\documentclass{article}
\usepackage[OT2,T1]{fontenc}
\DeclareSymbolFont{cyrletters}{OT2}{wncyr}{m}{n}
\DeclareMathSymbol{\Sha}{\mathalpha}{cyrletters}{"58}

\usepackage{amsbsy,amssymb,amscd,amsfonts,latexsym,amstext,delarray,
amsmath, diagrams}  \headsep=15pt
\usepackage[all,knot]{xy}
\def\cP{ {\cal P} }
\def\Un{ \mbox{Un} }

\def\cX{{\cal X} }

\def\End{\mbox{End}}
\def\bX{ \bar{X}}

\def\dim{{ \mbox{dim} }}
\def\Spec{{ \mbox{Spec} }}

\def\Hom{{ \mbox{Hom} }}

\def\cA{ {\cal A} }

\def\Res{ \mbox{Res} }
\def\om{{ \omega  }}

\def\ra{{ \rightarrow }}

\def\a{{ \alpha }}
\def\b{{ \beta} }
\def\g{{ \gamma }}
\def\d{{ \delta }}

\def\e{{ \epsilon }}
\def\F{ {\mathbb F} }

\def\Un{ \mbox{Un} }
\def\Vect{ \mbox{Vect} }
\def\hra{{ \hookrightarrow }}

\def\C{{ \mathbb{C} }}

\def\bs{ \backslash}

\def\G{{ \Gamma }}
\def\Gal{{ \mbox{Gal} }}

\def\bQ{\bar{\Q}}

\def\cA{ {\cal A}}

\def\Z{{ \mathbb{Z}}}

\def\bq{\begin{quote}}
\def\eq{\end{quote}}
\def\Aut{ \mbox{Aut}}

\newtheorem{thm}{Theorem}[section]
\newtheorem{cor}[thm]{Corollary}
\newtheorem{lem}[thm]{Lemma}

\def\Q{\mathbb{Q}}

\def\invlim{\varprojlim}

\def\be{\begin{equation}}
\def\ee{\end{equation}}

\def\Res{\mbox{Res}}

\def\D{ \Delta}

\def\max{\mbox{max}}

\def\om{\omega}

\def\tS{ \tilde{S}}
\def\loc{\mbox{loc}}

\def\b{ \bar}

\def\L{\Lambda}

\def\cP{{\cal P}}

\def\a{\alpha}

\def\bd{\begin{diagram} }
\def\ed{ \end{diagram}}

\def\b{\beta}

\def\uphi{\underline{\phi}}

\def\bs{\begin{slide}}
\def\es{\end{slide} }

\def\upsi{\underline{\psi}}
\def\bd{ \begin{diagram}}
\def\ed{\end{diagram} }
\def\Sym{\mbox{Sym}}

\def\cF{{\cal F}}
\def\te{\tilde{e}}
\def\tB{\tilde{B}}
\def\Res{\mbox{Res}}
\def\Gm{\mathbb{G}_m}

\def\cP{{\cal P}}
\def\ord{\mbox{ord}}

\def\ad{\mbox{ad}}
\title{Selmer varieties for curves with CM Jacobians}
\author{John Coates and Minhyong Kim}
\begin{document}
\maketitle
\begin{abstract}
We study the Selmer variety associated to a canonical quotient of the $\Q_p$-pro-unipotent fundamental
group of a smooth projective curve of genus at least two defined over
$\Q$ whose Jacobian decomposes into a product of abelian varieties
 with complex multiplication. Elementary multi-variable Iwasawa theory is used to
prove dimension bounds, which, in turn, lead to a new proof of Diophantine finiteness over $\Q$ for such curves.
\end{abstract}

Let $X/\Q$ be a smooth proper  curve of genus $g\geq 2$ and $b\in X(\Q)$ a rational point. We assume
that $X$ has good reduction outside a finite set $S$ of primes and choose an odd prime
$p\notin S$.
In earlier papers (\cite{kim1}, \cite{kim2}, \cite{kim3}, \cite{kim4}, \cite{KT}), a $p$-adic {\em Selmer variety}
$$H^1_f(G, U)$$
was defined and studied, with the hope of applying its structure theory to the
Diophantine geometry of $X$. Here, $G=\Gal(\bQ/\Q)$, $$U=\pi_1^{\Q_p,un}(\bX,b)$$
is the $\Q_p$-pro-unipotent
\'etale fundamental group of $$\bX=X\times_{\Spec(\Q)}\Spec(\bQ),$$ and the subscript $f$ refers
to a collection of local `Selmer conditions,' carving out a moduli space of torsors
for $U$ on the \'etale topology of $\Spec(\Z[1/S, 1/p])$ that satisfy the condition
of being {\em crystalline} at $p$.

The Selmer variety is actually a pro-variety consisting of a projective system
$$\cdots \ra H^1_f(G, U_{n+1}) \ra H^1_f(G, U_{n+1}) \ra \cdots \ra H^1_f(G, U_{2}) \ra H^1_f(G, U_{1})$$
of varieties over $\Q_p$ associated to the descending central series filtration
$$U=U^1\supset U^2\supset U^n\supset U^{n+1}=[U,U^n]\supset \cdots$$
of $U$ and the corresponding system of quotients
$$U_n=U^{n+1}\backslash U$$
that starts  out with
$$U_1=V=T_pJ\otimes \Q_p,$$
the $\Q_p$-Tate module of the Jacobian $J$ of $X$.

As a natural extension of the  map
$$\bd X(\Q) &\rTo& J(\Q) &\rTo& H^1_f(G,V)\ed$$
visible in classical Kummer theory,
the Selmer variety is endowed with a system of {\em unipotent Albanese maps} emanating from the
points of $X$:
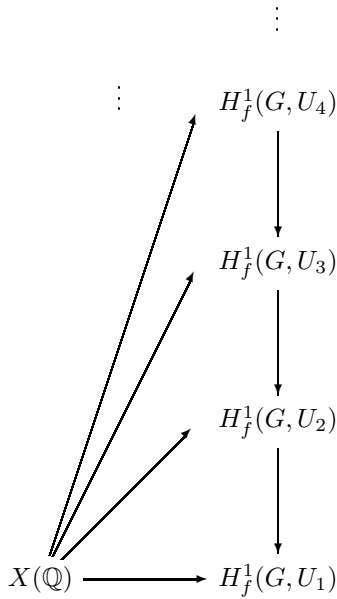
$$\begin{diagram}
 & &  & \vdots \\
  & \vdots &  &H^1_f(G,U_4) \\
  &\ruTo(2,6)&  &\dTo \\
 & & &H^1_f(G,U_3)\\
 &\ruTo(2,4) & &\dTo\\
 &  & &H^1_f(G,U_2) \\
 & \ruTo  & &\dTo\\
X(\Q) & \rTo& &H^1_f(G,U_1)
\end{diagram}$$
These maps  fit into commutative diagrams
$$\begin{diagram}
X(\Q) & \rTo & X(\Q_p)& & \\
\dTo & & \dTo & \rdTo & \\
H^1_f(G,U_n) & \rTo^{\loc_p}& H^1_f(G_p, U_n)& \rTo^{D} & U^{DR}_n/F^0
\end{diagram}$$
involving the local Selmer varieties $H^1_f(G_p, U_n)$ and their De Rham realizations $U^{DR}_n/F^0$. Here,
$F^0$ refers to the zeroth level of the Hodge filtration
$$U^{DR}\supset \cdots \supset F^i \supset F^{i+1} \supset \cdots \supset F^0$$
on the De Rham fundamental group $U^{DR}$ of $X\times_{\Spec(\Q)} \Spec(\Q_p)$.
Recall that the De Rham fundamental group is defined using the Tannakian
category of unipotent
vector bundles with flat connections on $ X\times_{\Spec(\Q)} \Spec(\Q_p)$ (\cite{kim2}, section 1),
and that the Hodge filtration $F^{-i}$ on $U^{DR}$ is the subvariety defined
by the ideal $F^{i+1}\cA^{DR}$ in the coordinate ring $\cA^{DR}$ of $U^{DR}$ (loc. cit. and \cite{deligne}). $F^0U^{DR}$ turns out to
 be a subgroup. The filtration on $\cA^{DR}$ is
defined over $\C$ using the $(d,\bar{d})$ decomposition on iterated integrals of differential forms,
but descends to any field of characteristic zero \cite{wojtkowiak}.
Here and in the following, we will suppress from the notation the object that
the Hodge filtration filters when the context provides sufficient clarity.

Diophantine
motivations oblige us to study the localization map $D\circ \loc_p$ with some care.
In fact, one could formulate the {\em dimension hypothesis}
$$\dim H^1_f(G, U_n) < \dim U^{DR}_n/F^0 \ \ \ \ \ \ \ \ \ \ \ (DH_n) $$
for each $n$, and show that $DH_n$ for any fixed $n$ implies the finiteness
of $X(\Q)$ \cite{kim2}.  Throughout this paper, dimension will refer to that of   algebraic
varieties over $\Q_p$ \cite{kim1}, although the dimensions of various associated graded
objects, e.g., $U^n/U^{n+1}$, are just the naive ones of $\Q_p$-vector spaces.
Given any $X$, it seems reasonable to believe that $DH_n$ should be true for
$n$ sufficiently large \cite{kim2}.

An eventual goal is to use the Selmer variety to arrive at a  structural understanding
of the Diophantine set $X(\Q)$, or at least some means of effective computation.
The hope for effective computation is associated with
classical method of Chabauty and Coleman \cite{coleman}, which the study of the unipotent Albanese map
generalizes. The related issue of structural understanding, on the other hand,
should concern an implication of the form
\bq
control of $L$-values $\Rightarrow$ control of Selmer varieties
\eq
following a pattern familiar from the theory of elliptic curves (\cite{CW}, \cite{rubin}).

It should be admitted right away that our current intuition for the
nature of such an implication is in a very tentative state. Nevertheless,
the cases studied previously  of hyperbolic curves of genus zero
and one  seem to suggest that our expectations are not entirely without ground.

The purpose of the present paper is to  augment our list of examples where something can be worked out
with
the case where $J$ is isogenous over $\bQ$ to a product
$$J\sim \prod_iA_i$$
of abelian varieties $A_i$
that have   complex multiplication by
CM fields $K_i$ of degree $2\dim A_i$. For this discussion we choose the prime $p$ to further
satisfy the condition that $p$ is split in the compositum $K$ of the fields $K_i$, and hence,
in each field $K_i$.

Let $\Q_T$ be the  maximal extension of $\Q$ unramified outside $T=S\cup \{p, \infty\}$
and $G_T=\Gal(\Q_T/\Q)$.
 Now let $$W=U/[U^2,U^2]$$ be the
quotient of $U$ by the third level of its derived series and let
$$W=U/[ U^2, U^2].$$
Of course $W$ itself has a descending central series
$$W=W^1\supset W^2\supset  \cdots \supset W^{n+1}=[W,W^n]\supset \cdots$$
and associated quotients $W_n=W/W^{n+1}$.
\begin{thm}
There is a constant $B$ (depending on $X$ and $T$) such that
$$\dim \sum_{i=1}^{n}H^2(G_T,W^i/W^{i+1})\leq B n^{2g-1}.$$
\end{thm}
We will derive this inequality as
 a rather elementary consequence of multi-variable Iwasawa theory.
 The key point is to control the distribution of zeros of a {\em reduced
 algebraic $p$-adic $L$-function} of sorts, namely an annihilator of
 a natural ideal class group.

In accordance with the `motivic nature' of the
construction, $W$ also has a De Rham realization $$W^{DR}=U^{DR}/[(U^{DR})^2,(U^{DR})^2]$$ over $\Q_p$,
endowed with a Hodge filtration.
 The upper bound of theorem 0.1 combines with an easy linear independence argument
 for sufficiently many elements in $(W^{DR})^n/(W^{DR})^{n+1}$, yielding a lower bound for
 the De  Rham realization $$
 W^{DR}_n/F^0$$
 of its local Selmer variety. We obtain thereby the
easy but important corollary:
\begin{cor}
For $n$ sufficiently large, we have the bound
$$\dim H^1_f(G, W_n) < \dim W^{DR}_n/F^0.$$
\end{cor}
Of course, this implies:
\begin{cor}(Faltings' theorem, special case.)
$X(\Q)$ is finite.
\end{cor}

For   explicit examples where the hypothesis is satisfied, we have of course the Fermat curves
$$x^m+y^m=z^m$$
for  $m\geq 4 $ (\cite{schmidt}, VI, Satz 1.2, Satz 1.5), but also the
twisted Fermat curves
$$ax^m+by^m=cz^m$$
for $a,b,c, \in \Q\setminus \{0\}$, $m\geq 4$. One might hope (optimistically)
that the methods of this paper will eventually lead to some effective understanding
of these twists.
Some relatively recent examples of hyperelliptic curves with CM Jacobians can be
found in
\cite{wamelen}. One from the list there is
$$y^2=-243x^6+2223x^5-1566x^4-19012x^3+903x^2+19041x-5882$$
whose Jacobian has CM by
$$\Q(\sqrt{-13+3\sqrt{13}}).$$

The results here conclude the crude application of
Selmer varieties to finiteness over $\Q$ in situations where the controlling Galois group of the base is
 essentially abelian. It remains then to work out the appropriate
interaction between  non-commutative geometric fundamental
groups and the non-commutative Iwasawa theory of number fields.

Of course, as far as a refined study of defining ideals for the
image of $D\circ \loc_p$ is concerned, work of any serious nature has not yet commenced.
In this regard, we note that there is little need in this paper for specific information
about the annihilator that occurs in the proof of theorem 0.1.
However, it is our belief that structure theorems of the `Iwasawa main conjecture' type
will have an important role to play in eventual refinements of the theory.
\section{Preliminaries on complex multiplication}

Let $F/\Q$ be a finite extension with the property that the isogeny decomposition
$$J\sim \prod_i A_i$$
as well as
the complex multiplication on each $A_i$ are defined over $F$. We assume further that  $F\supset \Q(J[4p])$, so that
$F_{\infty}:=F(J[p^{\infty}])$ has Galois group $\G\simeq \Z_p^r$
over $F$. Denote by $G_{F,T}$, the Galois group
$\Gal(\Q_TF/F)$.

As a representation of $G_{F,T}$,
we have $$V:=T_pJ\otimes \Q_p\simeq \oplus_{i} V_i$$
where $$V_i:=T_pA_i\otimes \Q_p.$$
Let $m$ be a modulus of $F$ that is divisible by the conductor of all the
representations $V_i$.
Each factor representation
$$\rho_i:G_{F,T}\ra (K_i\otimes \Q_p)^*\subset \Aut(V_i)$$
corresponds to an algebraic map
$$f_i:S_m\ra \Res^{K_i}_{\Q} (\Gm),$$
where $S_m$ is the Serre group of $F$ with modulus $m$ (\cite{serre2}, II) and
$\Res^{K_i}_{\Q}$ is the restriction of scalars from $K_i$ to $\Q$.
That is, there is a universal representation (op. cit. II.2.3)
$$\e_p:G_{F,T}\ra S_m(\Q_p)$$
such that
$$\bd \rho_i=f_i\circ \e_p:G_{F,T}& \rTo &S_m(\Q_p) &\rTo& \Res^{K_i}_{\Q} (\Gm)(\Q_p)=(K_i\otimes \Q_p)^*.\ed$$
Since we have chosen $p$ to split in each $K_i$, we have
$$\Res^{K_i}_{\Q} (\Gm)\otimes \Q_p \simeq \prod [\Gm]_{\Q_p}.$$
Each of the algebraic characters
$$\bd f_{ij}=pr_j\circ \rho_i: [S_m]_{\Q_p}& \rTo&  [\Res^{K_i}_{\Q_p} (\Gm)]_{\Q_p}\simeq \prod [\Gm]_{\Q_p} &\rTo^{pr_j}& [\Gm]_{\Q_p} \ed$$
correspond to Galois characters
$$\chi_{ij}=f_{ij}\circ \e_p: G_{F,T}\ra \Q_p^*$$
in such a way that
$$\rho_i\simeq \oplus_j \chi_{ij}.$$
Recall that $S_m$ fits into an exact sequence
$$0\ra T_m \ra S_m \ra C_m\ra 0$$
with $C_m$  finite and $T_m$ an algebraic torus (op. cit. II.2.2). Hence, there is an integer $N$ such that the kernel and cokernel of the restriction
map on characters
$$X^*([S_m]_{\Q_p}) \ra X^*([T_m]_{\Q_p})$$
is killed by $N$. On the other hand, $X^*([T_m]_{ \Q_p})$ is a finitely generated torsion free abelian
group. Let $\{\b_1', \ldots, \b_d'\}$ be a basis for the subgroup of $X^*([T_m]_{ \Q_p})$ generated by
the restrictions $f_{ij}|[T_m]_{ \Q_p}$ as we run over  all $i$ and $j$. Then the set $\{(\b_1')^N, \ldots, (\b_d')^N\}$
can be lifted to characters $\{\b_1, \ldots, \b_d\}$ of $[S_m]_{ \Q_p}$ so that
each $f_{ij}^{N^2}$ is a product
$$f_{ij}^{N^2}=\prod_k \b_k^{n_{ijk}}$$
for integers $n_{ijk}$. For ease of notation, we will now change the indexing
and write $\{f_1, \ldots, f_{2g}\}$ for the set of $f_{ij}$ and $\{\chi_i\}_{i=1}^{2g}$ for
the characters of $G_{F,T} $ that they induce.
 We have shown that
there are integers $n_{ij}$ such that
$$f_i^{N^2}=\prod_j \b_j^{n_{ij}}.$$
Thus, if we denote by $\xi_i$ the character
$$\b_i\circ \e_p: G_{F,T} \ra \Q_p^*,$$
then
$$\chi_i^{N^2}=\prod_j \xi_j^{n_{ij}}.$$
\begin{lem}
The characters $\xi_i$ are $\Z_p$-linearly independent.
\end{lem}
{\em Proof.}
The image of the map $\e_p: G_{F,T} \ra S_m(\Q_p) $  contains an open
subgroup $O_m$ of $T_m(\Q_p)$ (\cite{serre2} II.2.3, Remark).
 Suppose
$$\prod \xi_i^{a_i}=1$$
for some $a_i\in \Z_p$ as a function on $G_{F,T}$ (and, say, the choice of $p$-adic log
such that $\log(p)=0$).  Then
$$\prod_i \b_i^{a_i}=1$$
as a function on $O_m$.
Since the  $\b_i|[T_m]_{\Q_p}=(\b_i')^N$ are $\Z$-linearly independent, for each $j$, there
exists a cocharacter $$c_j: [\Gm]_{\Q_p} \ra [T_m]_{\Q_p}$$ such that
$\b_i\circ c_j=1$ for $i\neq j$ and
$$\b_j\circ c_j: [\Gm]_{\Q_p} \ra [\Gm]_{\Q_p}  $$ is non-trivial, and hence, an isogeny.
But
$c_j^{-1}(O_m)$ is an open subgroup of $\Q_p^*$. Hence, it contains an element of
the form $x=1+p^nu$, with $n>0$ and $u \in \Z_p^*$. Therefore,
$c=\b_j(c_j(x))\in 1+p\Z_p$ also has infinite order and
$c^{a_j}=1$.
Therefore, we get $a_j=0$.
$\Box$

Since the kernel of $$\rho=\oplus_j \rho_j=\oplus_{i=1}^{2g} \chi_i$$ is the same as that of
$\xi:=\oplus_{i=1}^d\xi_i$, $\xi$ maps $\G$ isomorphically
 to a subgroup of $\oplus_{i=1}^d(1+p\Z_p)$
of finite-index. After enlarging $F$ if necessary, we can assume that
there is a basis $\{\g_1, \ldots, \g_d\}$ for $\G$ such that
$\xi_i(\g_j)=1$ for $j\neq i$ and $\xi_i(\g_i)$ is a generator for
$\xi_i(G_{F,T})$. Here, we abuse notation a bit and write $\xi_i$ for the character of
$G_{F,T}$ as well as that of the quotient group $\G$ that it induces.

In the following, for any character $\phi$,  we will frequently
use the notation `$\phi $'  for the one-dimensional
vector space $\Q_p(\phi)$ on which $G_{F,T}$ acts via $\phi$,
as well as for the character itself.
Choose a basis $$B=\{e_1, e_2, \ldots, e_{2g}\}$$
of $V$ so that $e_i$ is a basis of $\Q_p(\chi_i)$. Write $\psi_i$ for the dual of
$\chi_i$.

Note that over $F$, the abelian variety $J$ has
good reduction everywhere \cite{ST}.
\section{Preliminaries on dimensions}

For a (pro-algebraic)  group or a Lie algebra $A$, we define the descending central series by
$$A^1=A; \ \ \ A^{n+1}=[A,A^n]$$
and the derived series by
$$A^{(1)}=A; \ \ \ A^{(n+1)}=[A,A^{(n)}].$$
The corresponding quotients are denoted by
$$A_n:=A/A^{n+1}$$
and
$$A_{(n)}:=A/A^{(n+1)}.$$
Also, we denote by $$Z_n(A):=A^n/A^{n+1}$$
the associated graded objects so that we have an exact sequence
$$0\ra Z_n(A) \ra A_n \ra A_{n-1}\ra 0.$$
Denote by
$$Z(A):=\sum_{n=1}^{\infty} Z_i(A)$$
the associated graded Lie algebra, described in \cite{serre}, II.1, in the case of a group.

According to \cite{quillen}, Appendix 3, the $\Q$-pro-unipotent completion of a finitely-presented discrete group $E$
can be constructed as follows: take the group algebra
$\Q[E]$, and complete it with respect to the augmentation ideal $K$:
$$\Q[[E]]:=\invlim_n \Q_p[E]/K^n.$$
Since the co-product
$$\D: \Q_p[E]\ra \Q_p[E]\otimes \Q_p[E]$$
defined by sending an element $g\in E$ to
$$g\otimes g\in \Q[E]\otimes \Q[E]$$ takes
$K$ to the ideal $$K\otimes \Q[E]+ \Q[E]\otimes K,$$
there is an induced co-product
$$\bd
\D: \Q[[E]]&\rTo& \Q[[E]]\hat{\otimes} \Q[[E]]:=\invlim_n (\Q[[E]]\otimes \Q[[E]])/(K\otimes \Q[E]+ \Q[E]\otimes K)^n.\ed$$
The unipotent completion
$U(E)$ can be realized as  the
group like elements in $\Q[[E]]$:
$$U(E)=\{g\in \Q[[E]] \ \ | \ \ \D(g)=g\otimes g\}.$$
This turns out to define the $\Q$-points of a pro-algebraic group over $\Q$.
Its Lie algebra,
$LieU(E)$, consists of the primitive elements
$$LieU(E)=\{X\in \Q[[E]] \ \ | \ \ \D(X)=X\otimes 1+1\otimes X\}.$$
For any element $g\in U(E)$,
$$\log(g)=(g-1)-(g-1)^2/2+(g-1)^3/3-\cdots$$
defines an element of $Lie U(E)$, and an elementary computation
starting from
$$\D(g-1)=(g-1)\otimes 1+1\otimes (g-1)+(g-1)\otimes (g-1)$$
shows that $\log (g) \in Lie U(E)$. In fact, this map is a bijection (loc. cit.)
$$\log: U(E) \simeq Lie U(E).$$

When $E$ is a topologically finitely presented pro-finite group, the $\Q_p$-pro-unipotent completion
$U_{\Q_p}(E)$ is defined in an entirely analogous manner, except
that the group algebra
$\Q_p[[E]]$ is defined somewhat differently:
First, let $E^{pro-p}$ be the  maximal pro-$p$ quotient of
$E$, and let $$\Z_p[[E^{pro-p}]]:=\invlim_{N} \Z_p[E^{pro-p}/N],$$ where
$N$ runs over the normal subgroups of $E^{pro p}$ of finite-index, be its Iwasawa algebra. Then
$$\Q_p[[E]]=\invlim_n [(\Z_p[[E^{pro-p}]]/K^n)\otimes \Q_p],$$
 where $K\subset \Z_p[[E^{pro-p}]]$ again denotes the augmentation ideal.
Then $$U_{\Q_p}(E)\subset \Q_p[[E]]$$ and its Lie algebra are defined exactly as above.
Consider the category $\Un(E,\Q_p)$ of unipotent continuous $\Q_p$-representations of
$E$, that is, finite-dimensional continuous representations
$$\rho:E \ra \Aut(M)$$
that possess a filtration
$$M=M^0\supset M^1\supset M^2\supset \cdots$$
such that each
$M^i/M^{i+1}$
is a direct sum of the copies of the trivial representation. We see that
$U$ acting on the left on $\Q_p[[E]]/K^n$
turns the system $\{\Q_p[[E]]/K^n\}$ into a pro-object of $\Un(E, \Q_p)$. Given any
pair $(M,m)$ where $M$ is a continuous unipotent $\Q_p$-representation of
$E$ and $m\in M$, there is a unique map of pro-representations
$$(\Q_p[[E]], e) \ra (M,m)$$
where $e \in \Q_p[[E]]$ comes from the identity of $E$, making
 the pair $(\Q_p[[E]]/K^n, e)$ universal among such pairs.  Therefore, if we let
 $$\bd F:\Un(E,\Q_p) & \rTo& \Vect_{\Q_p}\ed$$
be the forgetful functor from the category of unipotent continuous
$\Q_p$-representations of $E$ to the category of finite-dimensional $\Q_p$-vectors spaces,
the map
$$f \mapsto fe\in \Q_p[[E]]$$
defines an isomorphism
$$\End(F)\simeq \Q_p[[E]].$$
 Meanwhile, the condition of being group-like
corresponds to the compatibility with tensor products \cite{DM},
so that we have
$$U_{\Q_p}(E)=\Aut^{\otimes}(F),$$
the tensor-compatible automorphisms of $F$.

Since it will be our main object of interest,
we denote simply by $U$ the $\Q_p$-pro-unipotent completion of the pro-finite fundamental group
$\pi_1^{et}(\bX,b)$ of $\bX$  with  base-point at $b$.
Fix a  rational tangent vector $v\in T_bX$, and let $X'=X\setminus\{b\}$.
Let
$$U':=\pi_1^{\Q_p,un}(\bX',v),$$
the $\Q_p$-pro-unipotent completion of the profinite fundamental group of $\bX'$  with tangential base-point at $v$
as defined in \cite{deligne}.
These groups come with corresponding Lie algebras $L'=:Lie U'$ and $L:=Lie U$.

Let $$\Un(\bX,\Q_p)$$ be the category of unipotent lisse sheaves on the \'etale site of
$\bX$. Then  the fiber functor
$$F_b:\Un(\bX,\Q_p) \ra \Vect_{\Q_p},$$
which associates to any sheaf $\cF$ its stalk $\cF_b$, factors through the tensor equivalence of
categories
$$\Un(\bX, \Q_p) \simeq \Un(\pi^{et}_1(\bX,b), \Q_p)\stackrel{F}{\ra} \Vect_{\Q_p},$$
so that we also have
$$U=\Aut^{\otimes}(F_b).$$
Similarly,
$$U'=\Aut^{\otimes}(F'_b),$$
where
$$F'_b:\Un(\bX',\Q_p) \ra \Vect_{\Q_p}$$
is again the fiber functor defined by stalks.

As explained in \cite{hainmatsumoto}, Appendix A,
there are natural isomorphisms
$$U'\simeq U'_B\otimes \Q_p$$
and
$$U\simeq U_B\otimes \Q_p,$$
where $U'_B$ and $U_B$ denote the $\Q$-unipotent completions of
the topological fundamental groups $\pi'=\pi_1(X(\C), v)$ and $\pi=\pi_1(X(\C), b)$ of $X'(\C)$ and $X(\C)$.
(Appendix A of op. cit. is recommended in general for background on pro-unipotent
completions, while its section 2 contains a nice discussion of pro-algebraic groups.)
Therefore, $L'_B=:Lie U'_B$, and $L_B:=Lie U_B$ also satisfy  comparison isomorphisms
$$L'\simeq L'_B\otimes \Q_p$$
and
$$L\simeq L_B\otimes \Q_p.$$
The natural maps
$$\bd \pi' &\rTo& U'_B\ed$$
and
$$\bd \pi &\rTo& U_B\ed$$
induce isomorphisms
$$  Z(\pi')\otimes \Q \simeq Z(U'_B)$$
and
$$ Z(\pi)\otimes \Q \simeq Z(U_B).$$
(See, for example, \cite{amoros}, Prop. 1.2. In that reference, real coefficients are used. But this
obviously implies that same result for $\Q$-coefficients.)
On the other hand,
the bijections
$$\log: U'_B \ra L'_B$$
and
$$\log: U_B \ra L_B$$
take
$ghg^{-1}h^{-1}$ to
$[\log g, \log h]$
modulo higher commutators by the Baker-Campbell-Hausdorff formula (\cite{serre}, Chap. 4). Hence, there is a bijection of the descending central series filtrations on the two sides, and the brackets agree modulo terms of higher order.
Therefore, the log map also induces isomorphisms
$$Z(U'_B) \simeq Z(L'_B)$$
and
$$Z(U_B)\simeq Z(L'_B),$$
this time respecting the brackets.
From this, we also get
$$Z(\pi')\otimes \Q_p \simeq Z(L')$$
and
$$ Z(\pi)\otimes \Q_p \simeq Z(L).$$

It follows from this that $Z(L')$ is the free Lie algebra on $2g$ generators.
Hence, $L'$ is  free on generators obtained from any lift of a basis
for $(L')_1=V$. That is, we can take a lift $\tS$ of any basis $S$ for $V$, then the map
$$F(\tS)\ra L'$$
from the free Lie algebra on $\tS$ to $L'$ induces  isomorphisms
$$ F(\tS)/F(\tS)^n \simeq L'/(L')^n$$
 for each $n$, and hence, an isomorphism
 $$\overline{F(\tS)}:=\{F(\tS)/F(\tS)^n\}_n \simeq L'$$
 of pro-Lie algebras. (See, for example, \cite{hainmatsumoto}, appendix A again.)
 As generators for $L'$, we take a lifting
$\tB=\{\te_1, \cdots, \te_{2g}\}$ of the basis $B$ above.
The corresponding isomorphism
from  $\overline{F(\tB)}$ to $L'$  puts on $L'$ the structure of
a  graded pro-Lie algebra
$$L'=\oplus_{n=1}^{\infty} L'(n)$$
in such a way that
$$(L')^n=\oplus_{i\geq n}^{\infty} L'(i).$$
We warn the reader that this grading is {\em not} compatible
with the Galois action.
Since there appears to be little danger of confusion, we will denote the
elements $\te_i$ by $e_i$ again and the generating set $\tB$
by $B$.

  By \cite{labute}, the natural map
$$\pi' \ra \pi$$
induces an isomorphism
$$ (Z(\pi')\otimes \Q) /\bar{R}_B \simeq Z(\pi)\otimes \Q$$
for a Lie ideal $\bar{R}_B$ generated by the class of a single element $c:=\prod_i[a_i,b_i]\in (\pi')^2$, expressed in terms of a set of free generators
$\{a_1, \ldots,a_g, b_1, \ldots,b_g\}$
for  $\pi'$.
Consider the natural map
$$p:L' \ra L.$$
We have $\om:=\log (c)\in Ker(p)$, and the preceding discussion implies that
$$L'/R\simeq L,$$
where $R$ is the closed ideal generated by $\om$,
since there is induced an isomorphism of
associated graded algebras.

For the structure of $N:=LieW$, we  have therefore
$$N\simeq L'/[I+R],$$
where $$I=(L')^{(3)}=[[L',L'],[L',L']].$$

According to \cite{reutenauer}, we can construct a Hall basis for $L'$
as follows. First, we order $B$ so that $e_i< e_j$ if $i<j$.
This is, by definition, the set $H_0$. Now define $H_{n+1}$ recursively
as the brackets of the form
$$[\ldots[h_1,h_2],h_3],\ldots ], h_k]$$
where $k \geq 2$, $h_i\in H_n$, and
$$h_1<h_2\geq h_3\geq \cdots h_k.$$
Now choose a total order on $H_{n+1}$.
Finally, put $H=\cup_i H_i$ and extend the order by the condition
$$h\in H_i, \ \ k\in H_j,\ \ i<j \ \ \Rightarrow h>k.$$
Symbolically,
$$H_0>H_1>H_2>\cdots.$$
In fact, it is shown that $\cup_{i\geq n} H_i$ is a Hall basis for the subalgebra
$$(L')^{(n+1)}.$$
In particular, it follows that the elements of $H_1$ are linearly independent from
$(L')^{(3)}$, which is generated by $\cup_{i\geq 2}H_i$.
Furthermore, the basis consists of monomials,  so that
$H(i):=H\cap L'(i)$ is a basis for $L'(i)$.
Define $H_n(i):=H_n\cap L'(i)$ so that $H(i)=\cup_n H_n(i)$.
We thus get a bigrading
$$L'=\oplus L'(i,n),$$
where
$L'(i,n)$ is the span of $H_n(i)$.

Denote by $N'$ the Lie algebra
$$(L')_{(2)}=L'/I.$$
Then
\begin{lem}
For $n\geq 2$, the set $H_1(n)$, consisting of Lie monomials of the form
$$[[\ldots [e_{i_1}e_{i_2}]e_{i_3}]\ldots ]e_{i_n}],$$
where $i_1<i_2\geq i_3 \geq \cdots \geq i_n$, is linearly independent from
$$(L')^{(n+1)}+I.$$

\end{lem}
{\em Proof.}
We have the bigradings
$$I=\oplus_{i=1}^{\infty}\oplus_{j\geq 2} L_(i,j)$$
$$(L')^{(n+1)}=\oplus_{i\geq n+1}\oplus_{j=1}^{\infty} L(i,j),$$
from which it is clear that $(L')^{(n+1)}+I$ is the sum of
$L'(i,j)$ where $(i,j)$ runs over the pairs such
that $j\geq 2$ or $i\geq n+1$. Thus, $H_1(n)$ is
linearly independent from it. $\Box$
\medskip

\begin{cor}
The image $[H_1(n)]$ of $H_1(n)$ in $N'_n$ is a basis for
$Z_n(N')$.
\end{cor}
\medskip

The elements of $H_1(n)$ for $n\geq 2$ can be counted by noting that there
are $\binom{2g}{2}$ possibilities for the bracket $[e_{i_1}, e_{i_2}]$,
while for each such bracket, the cardinality of the non-increasing
$(i_3, i_4, \ldots, i_n)$ with $i_3\leq i_2$ is
$$\binom{(n-2)+(i_2-1)}{i_2-1}=\binom{n-3+i_2}{i_2-1}.$$
So we find the following dimension formula:
\begin{cor} For $n\geq 2$,
$$\dim Z_n(N')=\sum_{i=1}^{2g}(i-1)\binom{n-3+i}{i-1}.$$
\end{cor}
{\em Proof.} This follows immediately from the previous discussion together
with the observation that for any index $i$, there are $i-1$ possibilities
for the bracket $[e_j,e_i]$ at the beginning of an element of $H_1(n)$.
$\Box$
\medskip

We would like to understand the dimension of $Z_n(N)$.
Although it would be elementary work out a precise formula, we need
just a reasonable estimate for our purposes. That is, we need to estimate the dimension
of
$$Z_n(N')/ [Z_n(N')\cap Im(R)],$$
where $Im(R)$ refers to the ideal in $N'_n$ generated by the image of $\om$ (which we will
again  denote by $\om$).

For an ordered collection of elements
$v=(x_1, x_2, \ldots, x_m)\in B^m$  and an element $y\in L'$,
define
$$\ad_v(y):=[[ \ldots[y, x_1],x_2], \ldots, x_m]$$
Note that if $y\in (L')^2$ and $v'$ is a re-ordering of $v$, then
$$\ad_v(y)-\ad_{v'}(y)\in I$$
Thus, for $y\in (L')^2$, we have
$$\ad_v(y)\equiv \ad_{\ord(v)}(y)\ \  \mod \ I$$
where $\ord(v)$ is the unique reordering of $v$ for which the components are non-increasing.
Hence, any $x\in R\cap L'(n)$ has an expression as a linear combination
$$x\equiv \sum_i c_i \ad_{v_i}(\om) \ \ \ \mod I,$$
where $v_i$ runs through elements of $B^{n-2}$ with non-increasing components.
The number of such $v_i$ is
$$\binom{n-2+2g-1}{2g-1}=\binom{n-3+2g}{2g-1},$$
which therefore gives an upper bound on the dimension of
$Im(R)\cap Z_n(N')$.
\begin{lem}For $n\geq 2$,
$$\dim Z_n(N)\geq (2g-2)\binom{n-3+2g}{2g-1}+\sum_{i=1}^{2g-1}(i-1)\binom{n-3+i}{i-1}$$
\end{lem}

\section{Proofs}

We refer to \cite{kim1} and \cite{kim2}, section 2, for  general background material on Selmer varieties.
Recall that
$$H^1_f(G,W_n) \subset H^1(G_T, W_n)$$
consists of the cohomology classes corresponding to $W_n$-torsors
that are unramified outside $T$ and crystalline at $p$.
So a bound for $H^1(G_T, W_n)$  will be a bound
for $H^1_f(G_T, W_n)$ as well.

We will use again the exact sequence
$$0\ra H^1(G_T, Z_n(W)) \ra H^1(G_T, W_n) \ra H^1(G_T, W_{n-1}) $$
as in  op. cit. and a bound for the dimension of
$H^1(G_T, Z_n(W)) $. There is, as usual, the Euler characteristic formula \cite{jannsen} that
reduces over $\Q$ to
$$ \dim H^0(G_T, Z_n(W))-\dim H^1(G_T,Z_n(W))+\dim H^2(G_T,Z_n(W))=\dim[Z_n(W)]^{+}-\dim[Z_n(W)]$$
$$=-\dim[Z_n(W)]^{-},$$
where the signs in the superscript refer to the $\pm 1$ eigenspaces for the action of
complex conjugation. Because $Z_n(W)$ has weight $n$, we see that the $H^0$-term is zero
for $n\geq 1$, from which we get
$$(EC)\ \ \ \ \ \  \ \ \ \ \dim H^1(G_T,Z_n(W))=\dim[Z_n(W)]^{-}+\dim H^2(G_T,Z_n(W)).$$
Note that if $T'\supset T$, then
$$\dim H^1(G_{T'}, W_n)\geq \dim H^1(G_T, W_n).$$
The Euler characteristic formula then shows that
$$\dim H^2(G_{T'}, W_n)\geq \dim H^2(G_T, W_n)$$
as well. Therefore, in our discussion of bounds, we may increase the
size of $T$ to include the primes that ramify in the field $F$. In particular,
we may assume that $F\subset \Q_T$ so that
$$G_{F,T}\subset G_T,$$
a subgroup of finite index.

\medskip

{\em Proof of theorem 0.1}
 Since there is a constant in the formula, we can  assume
$n\geq 3$. Furthermore, by the surjectivity of the corestriction map
$$H^2(G_{F,T}, Z_n(N)) \ra H^2(G_T, Z_n(N)),$$
we may concentrate on bounding $H^2(G_{F,T}, Z_n(N))$.
As in \cite{kim4}, we consider the localization sequence
$$0\ra \Sha^2(Z_n(N)) \ra H^2(G_{F,T}, Z_n(N)) \ra \oplus_{v|T} H^2(G_v, Z_n(N)),$$
where $G_v=\Gal(\bar{F_v}/F_v)$.
For the local terms, we have Tate duality
$$H^2(G_v, Z_n(N))\simeq (H^0(G_v, Z_n(N)^*(1)))^*.$$
For $v\nmid p$, since $J$ has  good reduction, the action of
 $G_v$  on $Z_n^*(N)(1))$ is unramified.
But then, for $n\geq 3$, $Z_n(N)^*$ has Frobenius weight $\geq 3$, while $\Q_p(1)$
has Frobenius weight -2. Therefore,
$$H^0(G_v, Z_n(N)^*(1)) =0.$$
For $v|p$, we use instead the fact (\cite{fontaine}, theorem 5.2) that
$$H^0(H_v, Z_n(N)^*(1))=\Hom_{MF(\phi)}(F^{nr}_v, D_{cris}(Z_n(N)^*(1))).$$
Here, $F_v^{nr}$ is the maximal absolutely unramified subextension of  $F_v$
and $D_{cris}(\cdot)=((\cdot)\otimes B_{cris})^{G_v}$ is Fontaine's
crystalline Dieudonn\'e functor  applied to crystalline $G_v$-representations, while
$MF(\phi)$ is the category of admissible filtered $\phi$-modules over $F_v^{nr}$ (op. cit., section 5.1).
Since
each character $\psi_i$ occurs inside $H^1_{et}(\bar{J}, \Q_p)$, we know that $D_{cr}(\psi_i)$ occurs inside the
crystalline cohomology $H^1_{cr}(J, \Q_p)$ (\cite{FM}). But then, if the residue field of $F_v$ is of degree $d$ over $\F_p$,
$\phi^d$ again has positive weights on $D_{cris}(Z_n^*(N)(1))$ (\cite{KM}).
Therefore, $H^0(G_v, Z_n^*(N)(1))=0$.
It follows that
$$H^2(G_v, Z_n(N))\simeq \Sha^2(Z_n(N))\simeq [\Sha^1(Z_n(N)^*(1))]^*$$
by Poitou-Tate duality (\cite{milne}, theorem 4.10), where
$\Sha^1(Z_n(N)^*(1))$ is defined by the exact sequence
$$0 \ra \Sha^1(Z_n(N)^*(1)) \ra H^1(G_{F,T},Z_n(N)^*(1)) \ra \oplus_{v|T} H^1(G_v, Z_n(N)^*(1)).$$
Now, the group $\G=\Gal(F_{\infty}/F)$ is the image of $G_{F,T}$
inside $\Aut(J[p^{\infty}])$, and $Z_n(N)^*(1)$, being a sum of tensor products of
the characters $\psi_i=\chi_i^*$ and $\Q_p(1)$, is a direct summand of  $(V^*)^{\otimes n}(1)$.
Hence, by Bogomolov's theorem (as in op. cit. Lemma 6.20, Lemma 6.21),
$$H^1(\G, Z_n(N)^*(1))=0.$$
Therefore, using the Hochschild-Serre sequence, we get
$$H^1(G_T, Z_n(N)^*(1))\subset \Hom_{\L}(\cX_T, Z_n(N)^*(1))$$
for $$\L:=\Z_p[[\G]]\simeq \Z_p[[T_1, T_2,\ldots, T_{d}]]$$
and $\cX_T=\Gal(K_T/F_{\infty})$, the Galois group of the maximal abelian
pro-p extension $K_T$ of $F_{\infty}$ unramified outside $T$.
Here, $T_i=\g_i-1$ for free generators $\g_i$ of $\G$ chosen as in section 1 so that
$\xi_i(\g_j)=1$
while $\xi_i(\g_i)$ is a generator for the image of $\xi_i(G_{F,T})$.
The condition of belonging to the kernel of the localization map
will, in any case, imply
$$\Sha^1(Z_n(N)^*(1))\subset \Hom_{\L}(M, Z_n(N)^*(1))=\Hom_{\L}(M(-1), Z_n(N)^*),$$
where $M=M'/(\Z_p-\mbox{torsion})$ for the Galois group $M'=\Gal(H'/F_{\infty})$
of the $p$-Hilbert class field $H'$ of $F_{\infty}$. (Of course, we could take
an even smaller Galois group.)
According to \cite{greenberg0}, $M'$, and hence $M$, is a torsion $\Lambda$-module.
(That reference states this in the case where $F_{\infty}$ is replaced by
the compositum of all $\Z_p$-extensions of $F$, but the proof  clearly applies to
any $\Z_p^r$-extension.)
According to a lemma of Greenberg (\cite{greenberg}, Lemma 2), there is a subgroup $P\subset \G$ such that $\G/P\simeq \Z_p$ and
$M$ is still finitely generated over $\Z_p[[P]]$. Consequently, as explained in op. cit., page 89,
if we choose a basis $\{\e_1, \ldots, \e_{d-1}\}$ for $P$ and complete it to a basis of $\G$ using
an element $\e_{d}$ that maps to a topological generator of $\Z_p$, then in the variables $S_i=\e_i-1$,
we can take the annihilator to have the form
$$f=b_0(S_1,\ldots, S_{d-1})+b_1(S_1,\ldots, S_{d-1})S_{d}+\cdots +b_{l-1}(S_1,\ldots, S_{d-1})S_{d}^{l-1}+S_{d}^l$$
for some power series $b_i$. Furthermore, by approximation, we can choose the $\e_i$ to be of the
form
$$\e_i=\g_1^{n_{i1}}\cdots \g_{d}^{n_{i,d}}$$
for integers $n_{ij}$. That is, if the original $P$ involved $p$-adic powers
$n_{ij}$, we can approximate them by integral $n'_{ij}$ that are $p$-adically close, defining
another subgroup $P'$.
If the $S\subset M$ is a generating set as a $\Z_p[[N]]$ module,
then we see that $\Z_p[[P]]S=\Z_p[[P']]S$ mod $pM$, and hence, that
$\Z_p[[P']]S=M$ as well.

We know that $Z_n(N)$ is generated by the image of $H_1(n)$.
So $Z_n(N)^*$ is a subspace of the $G_{F,T}$-representation
given as the direct sum of the one-dimensional representations
$$\psi_{i_1}\otimes \psi_{i_2}\otimes \psi_{i_3} \otimes \cdots \otimes \psi_{i_n},$$
where $i_i$ run over indices from $\{1, 2, \ldots, 2g\}$ such that
$$i_1<i_2\geq i_3\geq \cdots \geq i_n.$$
So this is of the form
$$\oplus_{i<2g} [\psi_{i}\otimes \psi_{2g}\otimes \Sym^{n-2}(V^*)]\oplus K_n$$
where  $$\dim K_n \leq \binom{2g}{2}\binom{n-2+2g-2}{2g-2}=O(n^{2g-2}).$$
Therefore, the representation
$$\oplus_{i=3}^{n}Z_i(N)^*$$
is of the form
$$\oplus_{i<2g}[\psi_i\otimes \psi_{2g}\otimes (\oplus_{i=1}^{n-2}\Sym^i(V))]+R_n$$
where $R_n$ has dimension $\leq A n^{2g-1}$ for some constant $A$.
Clearly, $\dim H^2(G_T, R_n)\leq A' n^{2g-1}$ for another constant $A'$.
So we need to find a good bound for
$$\Hom_{\L}(M(-1), \oplus_{i<2g}[\psi_i\otimes \psi_{2g}\otimes (\oplus_{i=1}^{n-2}\Sym^i(V))]).$$
We  use the multi-index notation
$$\upsi^{\a}=\psi_1^{\a_1}\psi_2^{\a_2} \cdots \psi_{2g}^{\a_{2g}}$$
for a multi-index $\a=(\a_1, \ldots, \a_{2g})\in \mathbb{N}^{2g}$.
The weight of the multi-index $\a$ is denoted
$|\a|:=\sum_i \a_i$ so that
$$\Sym^i(V)=\oplus_{|\a|=i} \upsi^{\a}.$$
If a component
$$\Hom_{\L}(M(-1), \psi_i\otimes \psi_{2g}\otimes \psi^{\a})=\Hom_{\L}(M(-1)\otimes \chi_i\otimes \chi_{2g}, \psi^{\a})$$
is non-zero, then we must have
$$ \psi^{\a}(f_i)=0.$$
where
$$f_i:=f(c_{i1}S_1+c_{i1}-1, \cdots, c_{i, d}S_{d}+c_{i,d}-1)$$
$$=b_0^i(S_1,\ldots, S_{d-1})+b^i_1(S_1,\ldots, S_{d-1})S_{d}+\cdots +b^i_{l-1}(S_1,\ldots, S_{d-1})S_{d}^{l-1}+c_{id}^lS_{d}^l,$$
for some power series $b_j^i$ and units $c_{ij}:=\psi_i (S_j+1) \psi_{2g}(S_j+1)$, is in the annihilator  of $M(-1)\otimes \chi_i\otimes \chi_{2g}$. We wish to estimate how many zeros
each $f_i$ can have on the set $\{\upsi^{\a}\ | \ \ |\a|\leq n-2\}$.

There are independent elements $\{ \phi_i\} $ in the $\Z-$lattice of characters generated by
 the $\xi_i$  such that
$\phi_i(\e_j)=1$ for $i\neq j$
and
$$\xi_i=\phi_1^{m_{i1}}\cdots \phi_{d}^{m_{i,d}}$$
for  a nonsingular   matrix  $(m_{ij})$ with entries  $m_{ij}\in (1/M')\Z$ for some fixed denominator
$M'\in \Z \cap \Z_p^*$. Therefore, by the discussion in section 1,
we have
$$\psi_i=\phi_1^{q_{i1}}\cdots \phi_{d}^{q_{i,d}}$$
for a $(2g)\times d$ integral matrix $D=(q_{ij})$ of rank $d$ having entries in $(1/M)\Z$ for some fixed integer $M$.
Given a multi-index $\a$, we  then have
$$\upsi^{\a}=\uphi^{\a D},$$
with $\a D$ denoting the matrix product. For $|\a|\leq n-2$, we find the bound $$|\a D|\leq (n-2) (2g)|D|,$$
where $|D|=\max \{ |q_{ij}|\}$. Now, for each multi-index $$\d=(\d_1, \ldots, \d_d) \in [(1/M)\Z]^d$$ such that
$$|\d|\leq (n-2)(2g)|D|,$$ we need to count  the cardinality of the set
$$L_{\d}=\{\a \in \mathbb{N}^{2g} \ | \ \ \d=\a D,  |\a| \leq n-2\}.$$
If we fix one $\a\in L_{\d}$,  the map
$\a'\mapsto \a'-\a$ will inject $L_{\d}$ into the set of $\mu=(\mu_1, \ldots, \mu_{2g}) \in \Z^{2g}$  such that
$\mu D=0$ and $\sup_i|\mu_i|\leq (n-2)$. The first condition defines a lattice inside
a Euclidean space of dimension $2g-d$ while the second condition defines a fixed compact convex body (independent of
$n$) inside this space
dilated by a factor of $n-2$. Thus, there is a constant $C$ depending on the convex body such that
$|L_{\d}|\leq C(n-2)^{2g-d}$.
Now we turn to the number of  $\d$ for which
$$\uphi^{\d}(f_i)=0$$
and $|\d|\leq (n-2)(2g) |D|$. The  coefficients power series $b_j^i$ depend only
on $(\d_1, \ldots, \d_{d-1})$, which runs over a set of cardinality
$$\sum_{i=1}^{2g(n-2)M|D|} \binom{i+d-2}{d-2}.$$
This is the number of lattice points in   $d-1$ dimensional space that are contained inside
a simplex with vertices at the origin and at the points $$(0,  \ldots,0,2g(n-2)M|D| , 0, \ldots, 0).$$
This number is clearly majorized by the number of lattice points inside the cube $$[0, 2g(n-2)M|D|]^{d-1},$$
that is,
 $$(2g(n-2)M|D|+1)^{d-1}.$$ For each such $(\d_1, \ldots, \d_{d-1})$, there are at most $l$ $d$-tuples
$\d=(\d_1, \ldots, \d_{d-1}, \d_d)$ such that $\uphi^{\d}(f)=0$. We conclude that the number of
$\d$ such that $|\d|\leq 2g(n-2)M|D|$ and $\d(f)=0$ is bounded by
$ l(2g(n-2)M|D|+1)^{d-1}$.
Therefore, the number of zeros of each $f_i$ on $\{\upsi^{\a}\ \ | \ \ |\a|\leq n-2\}$
is bounded by $$C(n-2)^{2g-d} l(2g(n-2)M|D|+1)^{d-1}\leq An^{2g-1}$$
for some constant $A$.
For each such zero $\a$, the dimension
of $$\Hom_{\L}(M\otimes \chi_i\otimes \chi_{2g},\upsi^{\a})$$
is bounded by the number of generators $m$ of $M$.
From this, we deduce the desired asympototics
$$\dim \Hom_{\L}(M(-1), \oplus_{i<2g}[\psi_i\otimes \psi_{2g}\otimes (\oplus_{i=1}^{n-2}\Sym^i(V))])\leq mAn^{2g-1}. $$
$\Box$
\medskip

{\em Proof of corollary 0.2.}
For the rough estimates relevant to this paper, we will find  useful  the elementary fact that
$(n+a)^b=n^b+O(n^{b-1})$ for any fixed constant $a$ and exponent $b$.

We need to find lower bounds for the dimension of the local Selmer variety.
We have the De Rham realization
$$W^{DR}:=U^{DR}/(U^{DR})^{(3)},$$
where $U^{DR}$ is the De Rham fundamental group of $X\otimes \Q_p$ with base point at $b$.
We denote by $(U')^{DR}$ the De Rham fundamental group of $X'$ with basepoint
at $v$ (\cite{deligne}, \cite{kim2}).
Since $$(U')^{DR}\otimes \C_p\simeq U'\otimes \C_p,$$
(\cite{olsson1}, \cite{olsson2})
we see that
$$(L')^{DR}=Lie (U')^{DR}$$
is also free, and we can estimate dimensions exactly as in section 1.
For example, as in Lemma 1.3,
$$\dim Z_n(W^{DR})\geq (2g-2) \binom{n-3+2g}{2g-1}$$
so that
$$\sum_{i=3}^{n}\dim Z_i(W^{DR})\geq \frac{2g-2}{(2g)!}(n-2)^{2g}.$$

We need to estimate the contribution of $F^0(Z_n(W^{DR}))$.
For this, we let $\{b_1, \ldots, b_g, \ldots, b_{2g}\}$
be a basis of $(L')^{DR}_1$ such that $\{b_1, \ldots, b_g\}$ is a basis for
$F^0(L')^{DR}_1$. This determines a basis $H^{DR}=\cup_n H^{DR}_n$ for $(L')^{DR}$ following the
recipe of section 1. There are also corresponding bases
for $L^{DR}$, $(N')^{DR}:=Lie[(U')^{DR}/[(U')^{DR}]^{(3)}]$,
and a generating set for $N^{DR}$, exactly as in the discussion of section 1.
The Hodge filtration on $$(L')^{DR}_1=[H^1_{DR}(X'\otimes \Q_p)]^*$$ is
of the form
$$(L')^{DR}_1=F^{-1}(L')^{DR}_1\supset F^0(L')^{DR}_1 \supset F^1(L')^{DR}_1=0.$$
Hence, for an element
$$[[\ldots [b_{i_1}, b_{i_2}],b_{i_3}], \ldots] b_{i_n}]$$
of $H^{DR}_1(n)$ to lie in $F^0[Z_n((L')^{DR})]$, all of the $b_i$ must be in $F^0(L')^{DR}_1$.
 Thus, the dimension of $F^0[Z_n((N')^{DR})]$, and hence, of $F^0[Z_n(W^{DR})]$
is at most
$$\binom{g}{2}\binom{n+g-3}{g-1}.$$
From this, we get the estimate
$$\dim_{i=3}^{n} F^0(Z_i(W^{DR})) \leq cn^g$$
for some constant $c$. Therefore, we see that
$$(*) \ \ \ \ \dim W^{DR}_{n}/F^0=\dim W^{DR}_2/F^0+\sum_{i=3}^{n} Z_i(W^{DR})\geq \frac{(2g-2)}{(2g)!}n^{2g}+O(n^{2g-1}).$$

Now we examine the dimension of the minus parts
$Z_n(W)^{-}.$
For this, it is convenient to carry out the Hall basis construction with
yet another generating set. We choose $B'=\{f_1, \ldots, f_g, \ldots, f_{2g}\}$
so that $\{f_1, \ldots, f_g\}$ and $\{ f_{g+1}, \ldots, f_{2g}\}$ consist of
the plus and minus 1 eigenvectors in $V$, respectively.
Clearly, $\dim Z_n(N')^{-}$ will majorize $\dim Z_n(W)^{-}$.
Furthermore, as discussed above, $Z_n(N')=S_n+R_n$
where $S_n$ is the span of
$$[\ldots [f_j,f_{2g}],f_{i_3}]\ldots, f_{i_n}]$$
for $j<2g$ and nondecreasing $(n-2)$-tuples $(i_3, \ldots, i_n)$, while
$\dim R_n=O(n^{2g-2})$.
Now, $[f_j,f_{2g}]$ is in the minus part for $j\leq g$ and in the plus part for
$j\geq g+1$, while the contribution of the $(n-2)$-tuple will be as
$\Sym^{n-2}(V)$.

That is, $$\dim S_n^{-}=g \dim \Sym^{n-2}(V)^{+}+(g-1)\Sym^{n-2}(V)^{-}.$$

But
$$\Sym^{n-2}(V)=\oplus_i[\Sym^i(V^{+})\otimes \Sym^{n-2-i}(V^{-})]$$
of which we need to take into account the portions where $n-2-i$ is even and odd
respectively, to get the positive and negative eigenspaces.

For $n$ odd, we easily see that the plus and minus parts pair up, giving us
$$\dim \Sym^{n-2}(V)^{-}=\dim \Sym^{n-2}(V)^{+}=(1/2)\dim \Sym^{n-2}(V)=(1/2)\binom{n-3+2g}{2g-1}.$$
From this, we deduce that for $n$ odd,
$$\dim S_n^{-}=(1/2)(2g-1)\binom{n-3+2g}{2g-1}.$$
On the other hand, if $n$ is even, then
there is the embedding
$$\Sym^{n-2}(V)\hra \Sym^{n-1}(V)$$
$$v \ra v\cdot f_1$$
that preserves the plus and minus eigenspaces.
Hence,
$$\dim \Sym^{n-2}(V)^{+}\leq \dim \Sym^{n-1}(V)^{+}=(1/2)\binom{n-2+2g}{2g-1}=(1/2)n^{2g-1}/(2g-1)!+O(n^{2g-2})$$
and
$$\dim \Sym^{n-2}(V)^{-}\leq \dim \Sym^{n-1}(V)^{-}= (1/2)\binom{n-2+2g}{2g-1}=(1/2)n^{2g-1}/(2g-1)!+O(n^{2g-2}).$$
Therefore, for any $n$, we have
$$\dim S_n^{-}\leq (1/2)(2g-1)n^{2g-1}/(2g-1)!+O(n^{2g-2}).$$
and
$$\dim Z_n(N)^{-}\leq (1/2)(2g-1)\frac{n^{2g-1}}{(2g-1)!}+O(n^{2g-2}).$$
We deduce immediately that
$$\sum_{i=1}^{n} Z_i(N)^{-} \leq (1/2)(2g-1)\frac{n^{2g}}{(2g)!}+O(n^{2g-1}). $$
Combining this inequality with
the lower bound (*), theorem (0.1), and the Euler characteristic formula
$(EC)$, we get
$$\dim H^1_f(G,W_n)<\dim H^1_f(G_p, W_n)$$
for $n$ sufficiently large.
$\Box$
\medskip

{\bf Remark:} Note that in the comparison of leading coefficients,
$$(1/2)\frac{2g-1}{(2g)!} < \frac{2g-2}{(2g)!}$$
exactly for $g\geq 2$.
\medskip

{\em Proof of corollary 0.3}

By \cite{kim2}, section 2, and \cite{faltings}, there is an algebraic map
$$\bd D=D_{cr}: H^1_f(G_p, U)& \rTo& U^{DR}/F^0 \ed$$
sending a $U$-torsor
$$P=\Spec(\cP)$$
 to
$$\Spec(D_{cr}(\cP))=\Spec(\cP\otimes B_{cr})^{G_p},$$
an admissible $U^{DR}$ torsor, that is, a $U^{DR}$-torsor  with a compatible Frobenius action and a reduction of
structure group to $F^0U^{DR}$ (\cite{kim2}, section 1).

We wish to deduce an analogous map for $W$.
But \cite{olsson1} and \cite{olsson2} give
an isomorphism
$$L\otimes B_{cr}\simeq L^{DR}\otimes B_{cr}$$
compatible with the Lie algebra structure as well as the usual
Galois action, $\phi$-action, and Hodge filtration. In particular,
$$L^{(3)}\otimes B_{cr}\simeq (L^{DR})^{(3)}\otimes B_{cr},$$
and hence,
$$N\otimes B_{cr}\simeq N^{DR}\otimes B_{cr}.$$
Therefore,
$$D_{cr}(N)=N^{DR}$$
and
$$D_{cr}(W)=W^{DR}.$$
There is thereby
an induced map
$$\bd D:H^1_{f}(G_p, W) &\rTo & W^{DR}/F^0 \ed$$
following verbatim the construction for $U$ and $U^{DR}$
as in \cite{kim2}, section 2. That is, as in \cite{kim2}, Proposition 1 of section 1,
$W^{DR}/F^0$ classifies admissible torsors for $W^{DR}$, and the map assigns to
a $W$-torsor a $W^{DR}$-torsor, exactly following the recipe for $U$ and $U^{DR}$.

Now corollary (0.3) also follows verbatim the argument in \cite{kim1}, section 2,  and \cite{kim2}, section 3, by using the diagram
$$\begin{diagram}
X(\Q) & \rTo & X(\Q_p) & &\\
\dTo & & \dTo& \rdTo \\
H^1_f(G,W_n) & \rTo &H^1_f(G_p, W_n) & \rTo& W^{DR}_n/F^0$$
\ed.$$
for $n$ sufficiently large. We need only note that the map
$$]y[ \ra W^{DR}_n/F^0$$ from any residue disk $]y[ \subset X(\Q_p)$
to $W^{DR}_n/F^0$ has Zariski dense image, since the same
is true of
$$]y[ \ra U^{DR}_n/F^0$$
and the map
$$U^{DR}_n/F^0 \ra  W^{DR}_n/F^0$$
is surjective. $\Box$

\medskip

 {\bf Acknowledgements:} We are very grateful to Henri Darmon, Ralph Greenberg, Richard Hain,
Ramadorai Sujatha, and especially Kazuya Kato for many conversations related to the present subject.
 This paper was finalized during a stay at the Tata Institute for Fundamental Research
 in September, 2008. The hospitality of the staff and the organizers of the special semester
 on arithmetic geometry, E. Ghate, C.S. Rajan,
 and R. Sujatha, was essential to the completion of the present work.

{\footnotesize J.C.: Department of Pure Mathematics and Mathematical Statistics,
Centre for Mathematical Sciences,
University of Cambridge,
Wilberforce Road,
Cambridge,
CB3 0WB}

{\footnotesize M.K: Department of Mathematics, University College London,
Gower Street, London, WC1E 6BT, United Kingdom and The Korea Institute
for Advanced Study, Hoegiro 87, Dongdaemun-gu, Seoul 130-722, Korea}

\end{document}